\documentclass[12pt]{amsart}

\usepackage{amssymb}	% package for mathbb font
\usepackage{epsfig}		% package to include figures as eps files
\usepackage{amsthm}	% package to use \proof ... \qed syntax

\newcommand{\lcr}{\raisebox{-5pt}{\mbox{}\hspace{1pt}
                  \includegraphics{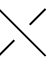}\hspace{1pt}\mbox{}}}
\newcommand{\ift}{\raisebox{-5pt}{\mbox{}\hspace{1pt}
                  \includegraphics{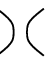}\hspace{1pt}\mbox{}}}
\newcommand{\zer}{\raisebox{-5pt}{\mbox{}\hspace{1pt}
                  \includegraphics{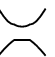}\hspace{1pt}\mbox{}}}

\newcommand{\torusL}{\raisebox{-7pt}{\mbox{}\hspace{-5pt}
                  \includegraphics{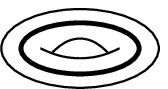}\hspace{-2pt}\mbox{}}}
\newcommand{\torusLM}{\raisebox{-7pt}{\mbox{}\hspace{-5pt}
                  \includegraphics{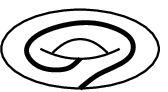}\hspace{-2pt}\mbox{}}}
\newcommand{\torusLMinv}{\raisebox{-7pt}{\mbox{}\hspace{-5pt}
                  \includegraphics{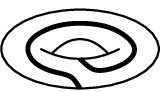}\hspace{-2pt}\mbox{}}}

\newcommand{\ot}{\otimes}
\newcommand{\tr}{\mathrm{tr}}
\newcommand{\bbc}{\mathbb{C}}
\newcommand{\slc}{SL_2(\mathbb{C})}
\newcommand{\LittleOneHalf}{{\textstyle \frac{1}{2}}}

\newtheorem{theorem}{Theorem}
\newtheorem{lemma}{Lemma}
\newtheorem{example}{Example}
\newtheorem{definition}{Definition}
\newtheorem{proposition}{Proposition}
\newtheorem{conjecture}{Conjecture}
\newtheorem*{remark}{Remark}

\title[Detecting Torsion in Skein Modules]{Detecting Torsion in Skein Modules using Hochschild Homology}

\author{Michael McLendon}
\address{Department of Mathematics and Computer Science,
Washington College, Chestertown, Maryland, 21620}
\email{mmclendon2@washcoll.edu}

\begin{document}

\begin{abstract}
Given a Heegaard splitting of a closed $3$-manifold, the
skein modules of the two handlebodies are modules over the skein
algebra of their common boundary surface.
The zeroth Hochschild homology of the skein
algebra of a surface with coefficients in the tensor product
of the skein modules of two handlebodies is interpreted
as the skein module of the $3$-manifold obtained by gluing the
two handlebodies together along this surface.  A
spectral sequence associated to the Hochschild complex
is constructed and conditions are given for the existence
of algebraic torsion in the skein module of this $3$-manifold.
\end{abstract}

\maketitle

\section{Introduction}
Skein modules were introduced independently by Przytycki \cite{Przytycki1991}
and Turaev \cite{Turaev1988}
and have been an active topic of research since their introduction.
In particular, skein modules underlie quantum
invariants \cite{Lickorish1993, KauffmanLins} and are connected to the representation
theory of the fundamental group of the manifold \cite{Bullock1997, PrzytyckiSikora2000}.

The skein module is spanned by the equivalence classes of framed links
in the $3$-manifold.  The skein module
of the cylinder over a surface has a multiplication that
comes from laying one framed link on top of the other.  With this
multiplication, the skein module of the cylinder over a surface
is an algebra.

Given a Heegaard splitting $H_0 \cup_{F} H_1$
of a $3$-manifold, the
skein module of each handlebody $B_i = K(H_i)$
is a module over the skein algebra
of the cylinder over the common boundary surface $A = K(F)$.
If $H_0$ is glued to $F \times \{0\}$ and $H_1$ is
glued to $F \times \{1\}$, then $B_0$ is a left $A$-module
and $B_1$ is a right $A$-module.  Hence
$B_0 \ot B_1$ is a bimodule over $A$.

We will use the interplay between the two
handlebodies and their common boundary surface to compute
the Hochschild homology of $A$ with coefficients in $B_0 \ot B_1$.
We will then construct a spectral sequence and show how it can
be used to detect algebraic torsion in the skein module
of the manifold.

\section{Preliminaries}

\subsection{Skein Modules}
Let $\mathcal{L}(M)$ denote the set of isotopy classes of framed links in $M$,
including the empty link, $\phi$.
Let $R=\mathbb{Z}[t,t^{-1}]$ be the ring of Laurent polynomials.
Consider the free module
$R \mathcal{L}(M)$ with basis $\mathcal{L(M)}$.  Define $S(M)$ to
be the smallest submodule of $R \mathcal{L}(M)$ containing
all expressions of the form
$\displaystyle{\lcr-t\zer-t^{-1}\ift}$
and
$L \sqcup \bigcirc + (t^2 + t^{-2}) L$, where $L$ is any framed link
and where the framed links in each expression are identical outside
the regions shown in the diagrams.  The
\emph{Kauffman bracket skein module}
$K(M)$ is the quotient $R \mathcal{L}(M)/S(M).$

Because $K(M)$ is defined using local relations on framed
links, two homeomorphic manifolds have
isomorphic skein modules.  Thus $K(M)$ is an invariant of the
$3$-manifold $M$.

\subsection{Heegaard Splittings}
Let $M$ be a closed, orientable, connected 3-manifold.
Then for some non-negative integer $g$ there exist
genus $g$ handlebodies $H_0$ and $H_1$ such that
$H_0 \cap H_1 = \partial H_0 = \partial H_1 = F$
is a closed, orientable, connected genus $g$ surface
and $H_0 \cup H_1 = M$.
We call these two handlebodies a
\textit{Heegaard splitting} of the manifold.
A simple proof of the existence of Heegaard splittings
using a triangulation of the manifold
can be found in Rolfsen \cite{Rolfsen}.
A given manifold may have many different
Heegaard splittings.

Note that we can take a neighborhood
of the surface $F$ and think of the Heegaard splitting as breaking
the manifold into $H_0$, $F \times [0,1]$, and $H_1$, where $H_0$ is
glued to $F \times \{0\}$ by the identity map
and $H_1$ is glued to $F \times \{1\}$ by a gluing map $f$.
We will model Heegaard splittings in this way, and we will be
interested in the properties of the gluing map $f$.

\subsection{Hochschild Homology}
Hochschild homology is a functor that associates an ordered
collection of $R$-modules to an $R$-algebra $A$ and an
$A$-bimodule $B$.  The Hochschild chain complex has chains
$C_n$ given by
\[
C_n = C_n(A; B) = B \ot (A^{\ot n}) = B \ot 
\underbrace{A \ot A \ot \dots \ot A}_{n~\mathrm{times}}
\]
for $n \geq 0$ and
$C_n = 0$ for $n < 0$.
The Hochschild boundary map $d_n : C_n \to C_{n-1}$ is given by
\begin{eqnarray}
d_n ( b \ot a_1 \ot \dots \ot a_n )
& = & b a_1 \ot a_2 \ot a_3 \ot \dots \ot a_n \nonumber \\
& & - b \ot a_1 a_2 \ot a_3 \ot \dots \ot a_n \nonumber \\
& & + b \ot a_1 \ot a_2 a_3 \ot \dots \ot a_n \nonumber \\
& & + \dots + (-1)^{n-1} b \ot a_1 \ot \dots \ot a_{n-1} a_n \nonumber \\
& & + (-1)^{n} a_n b \ot a_1 \ot a_2 \ot \dots \ot a_{n-1} \nonumber
\end{eqnarray}
where the products $a_i a_{i+1}$ take place in the algebra $A$
and the products $b a_1$ and $a_n b$ come from the respective
right and left actions on $B$ by $A$.
The Hochschild homology of $A$ with coefficients in $B$ is the homology
of the Hochschild complex and is denoted $HH_{*}(A; B)$.
If $B = A$ we will denote $HH_{*}(A; A)$ by $HH_{*}(A)$.

\section{The Hochschild Homology of a Heegaard Splitting}

\subsection{The skein module of a Heegaard splitting}
\begin{lemma}(Hoste-Przytycki \cite{HostePrzytycki1993})
Consider the manifold $F \times [0,1]$ where $F$ is a surface.
Let $\alpha$ be a simple closed curve on $F \times \{ 0 \}$ (or $F \times \{ 1 \}$).
Let $H$ be the manifold obtained by attaching a $2$-handle to $\alpha$.
Let $I$ be the submodule of $K(F)$ generated by relations of the form
$\{ s - h(s) \}$ where $s \in K(F)$ is a skein and $h(s)$ is the skein $s$
modified by a handleslide across $\alpha$.  Then $K(H) = K(F)/I$.
\label{handleslide}
\end{lemma}

The lemma above generalizes to the case where the manifold is obtained by
attaching more than one $2$-handle to $F \times [0,1]$.  Say we attach
a $2$-handle to $F \times \{ 0 \}$ along $\alpha$ and attach another
$2$-handle to $F \times \{ 1 \}$ along $\beta$.  Call the resulting
manifold $M$.  If $I$ is the submodule
of $K(F)$ generated by handleslides along $\alpha$ and $J$ is the submodule
of $K(F)$ generated by handleslides along $\beta$, then
$K(M) = K(F)/(I+J)$.
We can use this result and a property of the tensor product to
get the following proposition.
\begin{proposition}  (discussed by Frohman-Gelca in \cite{FrohmanGelca2000})
Let $M$ be a closed, connected, oriented $3$-manifold with Heegaard splitting
$M= H_0 \cup H_1$, $F = H_0 \cap H_1$.  Then
$K(M) = K(H_1) \ot_{K(F)} K(H_0)$.
\label{kmhs}
\end{proposition}

\proof
$H_0$ is obtained from
$F \times [0,1]$ by attaching $2$-handles to $F \times \{0\}$
along attaching curves $\alpha_k$.  Likewise, $H_1$ is obtained
from $F \times \{1\}$ by attaching $2$-handles to along curves $\beta_n$.
For $i \in \{0,1\}$, let $S_i$ be the submodule of $K(F)$
generated by handleslides across the $\alpha_k$ or across the $\beta_n$, respectively.
We can apply Lemma \ref{handleslide} to each $H_i$ and to $M$.  Then
we have $K(H_i) = K(F)/S_i$ and $K(M) = K(F)/(S_1 + S_0)$.

We know that $A/I \ot_A B \cong B/(IB)$ from homological algebra,
see Osborne \cite[Proposition 2.2]{Osborne}.
Let $A = K(F)$.  Consider $A/S_1 \ot_A A/S_0 \cong \frac{A/S_0}{S_1 (A/S_0)}$.
An element from $S_1 (A/S_0)$ looks like $sa + S_0$ where $s \in S_1$
and $a \in A$.
An element of $\frac{A/S_0}{S_1 (A/S_0)}$ looks like
$(a^{\prime} + S_0) + (sa + S_0)$.
Recall that the empty skein $\phi$ is the multiplicative
identity in $A$, thus $sa$ runs over all of $S_1$ and so
$(a^{\prime} + S_0) + (sa + S_0) = a^{\prime} + (S_1 + S_0)$.
Thus $K(H_1) \ot_{K(F)} K(H_0) \cong A/S_1 \ot_A A/S_0
\cong \frac{A/S_0}{S_1 (A/S_0)} \cong A/(S_1 + S_0) \cong K(M)$.
\qed

\subsection{Connection with character varieties}

Another interesting and useful property of the skein module $K(M)$ comes
when we specialize at $t=-1$.  The coordinate ring of
the $SL_2(\mathbb{C})$-characters on $\pi_1(M)$
is a quotient of this specialization.
This approach has been 
developed by Bullock \cite{Bullock1997} and also by
Przytycki and Sikora \cite{PrzytyckiSikora2000}.

Denote the specialization of $K(M)$ at $t=-1$ by $K_{-1}(M)$
and denote the space of $SL_2(\mathbb{C})$-characters by $X(M)$.
Let $\mathbb{C}^{X(M)}$ denote the algebra of functions on
$X(M)$ and let $R(M)$ denote the coordinate ring of $X(M)$.

By a theorem of Culler and Shalen \cite{CullerShalen1983},
$X(M)$ is an affine algebraic set.  Hence one can consider
the ring of polynomial functions on $X(M)$.  This ring of polynomial
functions is called the \textit{coordinate ring} of $X(M)$.
Indeed, Culler and Shalen show that the coordinate ring is
finitely generated.
% see Prop 1.4.1, p. 116 and Cor. 1.4.5 , p. 119 of Culler and Shalen

An oriented knot in $M$ determines a conjugacy class in $\pi_1(M)$
and thus an oriented knot $L$ defines a function
$\varphi_L : X(M) \to \mathbb{C}$ by $\varphi_L (\chi_{\rho}) = \chi_{\rho}(L) = \tr(\rho(L))$
where $\chi_{\rho}$ is the character induced by the representation $\rho$
and the knot $L$ is seen as an element of $\pi_1(M)$.
Since $\tr(A) = \tr(A^{-1})$ for any matrix $A$, the particular orientation on
the knot $L$ is irrelevant.

Let $\mathbb{C}\mathcal{L}(M)$ denote the vector space of
framed links in $M$.  Define a function
$\tilde{\Phi} : \mathbb{C}\mathcal{L}(M) \to \mathbb{C}^{X(M)}$ by
$\tilde{\Phi}(L) = - \varphi_L$ for a knot $L$ and
$\tilde{\Phi}(L) = \prod_i (-\varphi_{L_i})$ for a link
$L$ with components $L_i$.

\begin{lemma}(Bullock \cite{Bullock1997})
The map $\tilde{\Phi}$ descends to a map
$\Phi : K_{-1}(M) \to \mathbb{C}^{X(M)}$.  Its
image is the coordinate ring $R(M) \subset \mathbb{C}^{X(M)}$
and its kernel is the nilradical of $K_{-1}(M)$.
\label{bullock}
\end{lemma}

The proof that the map descends follows from the observation that
the skein relation maps to the $SL_2(\mathbb{C})$ trace identity
$\mathrm{tr}(AB) + \mathrm{tr}(AB^{-1}) = \mathrm{tr}(A)\mathrm{tr}(B)$.

Przytycki and Sikora \cite{PrzytyckiSikora2000} have shown that
the nilradical of $K_{-1}(M)$ is trivial for surfaces and handlebodies.
Thus for surfaces and handlebodies $\Phi$ is an isomorphism between the specialized
skein module and the coordinate ring of the character variety.

\begin{example}
As an example, let's look at the $3$-manifold $M = T^2 \times [0,1]$.  We
know that 
$\pi_1(M) = \langle \ell, m~|~\ell m \ell^{-1} m^{-1} = 1 \rangle$.
The coordinate ring $R(M)$ is generated by
$x = -\mathrm{tr}(\rho(m))$, $y = - \mathrm{tr}(\rho(\ell))$,
and $z = -\mathrm{tr}(\rho(\ell m))$, and it has one
relation induced by $\mathrm{tr}(\ell m \ell^{-1} m^{-1}) = 2$.
Hence
\[
R(M) \cong \bbc[x,y,z] / (x^2 + y^2 + z^2 + xyz - 4).
\]
\label{ex-torus}
\end{example}

% because the above section ends with a formula, there is too much
% space between the sections.  add a little negative space to correct this
% \vspace{-2\li}

\subsection{Using Hochschild Homology}
\label{hh-heegaard}

Now we can use the connection between the specialized skein
module and the coordinate ring in the context of a Heegaard
splitting.

Let $M$ be a 3-manifold with Heegaard
splitting $M = H_0 \cup F \times [0,1] \cup H_1$
and with gluing maps $f_0 : H_0 \to F \times \{0\}$ and
$f_1: H_1 \to F \times \{1\}$.
We will take $f_0$ to be the identity, hence the structure of the manifold
is described by $f_1$.

\begin{figure}[ht]
  \begin{center}
    \leavevmode
    \epsfxsize = 8cm
    \epsfysize = 4.5cm
    \epsfbox{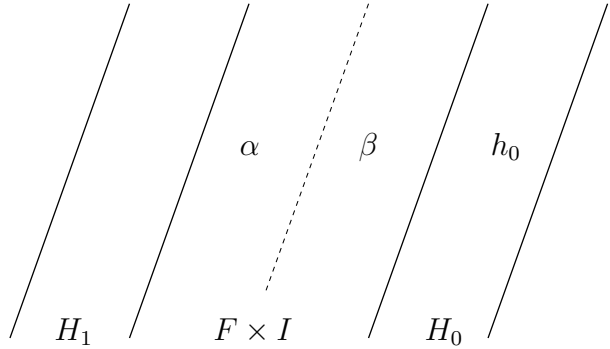}
    \put(-140,70){$\alpha$}
    \put(-95,70){$\beta$}
    \put(-45,70){$h_0$}
    \put(-150,0){$F \times I$}
    \put(-70,0){$H_0$}
    \put(-210,0){$H_1$}
    \caption{$(\alpha \beta) * h_0 = \alpha * (\beta * h_0)$ defines a left action}
    \label{left-action}
  \end{center}
\end{figure}

\begin{figure}[ht]
  \begin{center}
    \leavevmode
    \epsfxsize = 8cm
    \epsfysize = 4.5cm
    \epsfbox{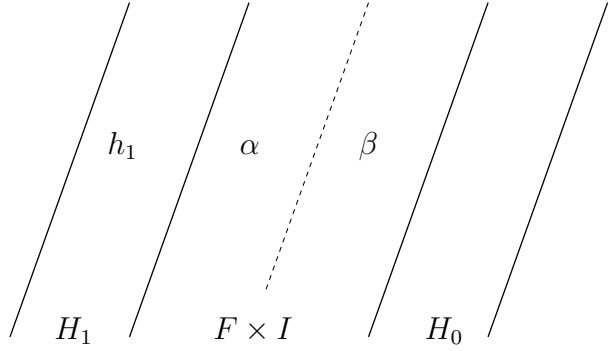}
    \put(-140,70){$\alpha$}
    \put(-95,70){$\beta$}
    \put(-190,70){$h_1$}
    \put(-150,0){$F \times I$}
    \put(-70,0){$H_0$}
    \put(-210,0){$H_1$}
    \caption{$h_1 * (\alpha \beta) = (h_1 * \alpha) * \beta$ defines a right action}
    \label{right-action}
  \end{center}
\end{figure}

The action of $K(F)$ on $K(H_i)$ is
given by pushing the skeins from $F \times I$ into $H_i$ using the inverse of
the $f_i$ gluing map.
The action of $\alpha, \beta \in K(F)$ on $h_0 \in K(H_0)$ is a left action,
$(\alpha \beta) * h_0 = \alpha * (\beta * h_0)$,
as shown in Figure \ref{left-action}.
The action of $\alpha, \beta \in K(F)$ on $h_1 \in K(H_1)$ is a right action,
$h_1 * (\alpha \beta) = (h_1 * \alpha) * \beta$, as shown in Figure
\ref{right-action}.
For $R = \mathbb{Z}[t, t^{-1}]$, we know that $A = K(F)$ is an algebra over
$R$.  Also, $B_0 = K(H_0)$,
and $B_1 = K(H_1)$ are modules over $R$.  Since $B_0$ is a
left $A$-module and $B_1$ is a right $A$-module, the tensor product
$B = B_0 \ot_R B_1$ is a bimodule over $A$.
We will use the unspecified tensor $\ot$ to denote $\ot_R$.

Now we look at $HH_{*}(A; B_0 \ot B_1)$.
With the choice of $B = B_0 \ot B_1$, the chains become
\[
C_n(A; B_0 \ot B_1) = (B_0 \ot B_1) \ot (A^{\ot n})
\]
for $n \geq 0$ and
$C_n = 0$ for $n < 0$.
Rearrange the terms in the tensor product
so that the $C_n$ become
\[
C_n(A; B_0 \ot B_1) = B_1 \ot
\underbrace{A \ot A \ot \dots \ot A}_{n~\mathrm{times}} \ot B_0.
\]
Then $d_n$ is
\begin{eqnarray}
d_n ( b_1 \ot a_1 \ot \dots \ot a_n \ot b_0 )
& = & b_1 a_1 \ot a_2 \ot a_3 \ot \dots \ot a_n \ot b_0 \nonumber \\
& & - b_1 \ot a_1 a_2 \ot a_3 \ot \dots \ot a_n \ot b_0 \label{bound-map} \\
& & + \dots + (-1)^{n} b_1 \ot a_1 \ot \dots \ot a_n b_0. \nonumber
\end{eqnarray}
Thus, $d_n : B_1 \ot (A^{\ot n}) \ot B_0
\to B_1 \ot (A^{\ot (n-1)} ) \ot B_0$.

Notice that without the $B_1$, the sequence
$$\dots \to A \ot A \ot A \ot B_0 \to A \ot A \ot B_0 \to
A \ot B_0 \to B_0 \to 0$$ is a free (hence projective) resolution of $B_0$.
If we delete $B_0$ from this complex, tensor over $A$ on the left by $B_1$
and compute the homology, we get $\mathrm{Tor}_i^A(B_1, B_0)$.
That is, $\mathrm{Tor}_i^A(B_1,B_0)$ is the homology of the following complex
$$\dots \to B_1 \ot_A A \ot A \ot A \ot B_0 \to B_1 \ot_A
A \ot A \ot B_0 \to B_1 \ot_A A \ot B_0 \to 0.$$
Since $B_1 \ot_A A = B_1$, the complex becomes
$$\dots \to B_1 \ot A \ot A \ot B_0 \to B_1 \ot A \ot B_0
\to B_1 \ot B_0 \to 0.$$
Therefore the $\mathrm{Tor}$ complex is exactly the same as the
Hochschild complex, that is,
\[HH_i(A; B_0 \ot B_1) = \mathrm{Tor}_i^A(B_1, B_0).\]

\begin{lemma}
The zeroth Hochschild homology of a Heegaard splitting
is the skein module of the 3-manifold $M$.
\end{lemma}

\proof
$\mathrm{Tor}_0$ corresponds to $\ot$, thus
$\mathrm{Tor}_0^A(B_1, B_0) = B_1 \ot_A B_0$.
We know from Proposition \ref{kmhs} that $K(M) = K(H_1) \ot_{K(F)} K(H_0)$, thus
$K(M) = B_1 \ot_A B_0 = \mathrm{Tor}_0^A(B_1, B_0) = HH_0(A; B_0 \ot B_1)$.
\qed

\section{A Spectral Sequence to Detect Torsion}
Next we use a filtration on $R$, $A$, $B_0$, and $B_1$ to get a
spectral sequence and search for $(1+t)$-torsion in $K(M)$.  We will follow
a process used by Brylinski in \cite{Brylinski1988} to study Poisson
manifolds.

The ring $R = \mathbb{Z}[t,t^{-1}]$ of Laurent polynomials has a
filtration by the ideals corresponding to powers of $(1+t)$.
\[\dots \subset (1+t)^3 R \subset (1+t)^2 R \subset (1+t) R \subset R\]
This is a decreasing filtration.  By manipulating the indices, we can
use this filtration to get an increasing filtration.  In particular, if
$\mathcal{F}_s(R) = (1+t)^{-s} R$ for $s \leq 0$ and
$\mathcal{F}_s(R) = R$ for $s > 0$, then $\mathcal{F}$ is
an increasing filtration on $R$.
This filtration extends to the $R$-modules $A=K(F)$ and $B_i=K(H_i)$ by
\[\dots \subset (1+t)^3 A \subset (1+t)^2 A \subset (1+t) A \subset A\] and
\[\dots \subset (1+t)^3 B_i \subset (1+t)^2 B_i \subset (1+t) B_i \subset B_i.\]
It also extends to the Hochschild complex $C_n = C_n(A;B_0 \ot B_1)$ by the
following from Brylinski \cite{Brylinski1988}.
\begin{eqnarray}
\mathcal{F}_s(C_n) & = & \mathcal{F}_s \left( B_1 \ot (A^{\ot n})
                         \ot B_0 \right) \nonumber \\
& = & \sum_{s_0 + \dots + s_{n+1} \leq s} \left( \mathcal{F}_{s_0}(B_1) \ot
      \mathcal{F}_{s_1}(A) \ot \dots \ot 
      \mathcal{F}_{s_n}(A) \ot \mathcal{F}_{s_{n+1}}(B_0) \right) \nonumber \\
& = & \sum_{\sum s_i \leq s} (1+t)^{-\sum s_i} \left( B_1 \ot (A^{\ot n})
      \ot B_0 \right) \nonumber
\end{eqnarray}

Now we create a spectral sequence $\{ E^r \}$ beginning at the
$E^0$ level with
\[E^0_{p,q} = \mathcal{F}_{p}(C_{p+q}) / \mathcal{F}_{p-1}(C_{p+q})\]
and the $E^0$ level boundaries $\Delta^0_n$ are induced by the Hochschild boundary
map $d_n$ as described in Equation \ref{bound-map}.

We move from the $E^0$ level to the $E^1$ level by taking
homology and letting the new boundary map $\Delta^1_n$ be the
connecting homomorphism induced by the short exact sequence
\[
0 \to \frac{\mathcal{F}_{p-1}(C_n)}{\mathcal{F}_{p-2}(C_n)} \to
\frac{\mathcal{F}_{p}(C_n)}{\mathcal{F}_{p-2}(C_n)} \to
\frac{\mathcal{F}_{p}(C_n)}{\mathcal{F}_{p-1}(C_n)} \to 0.
\]
In general,
 $E^r_{p,q} = H(E^{r-1}_{p,q}) = \mathrm{ker}(\Delta^{r-1})/\mathrm{im}(\Delta^{r-1})$ and
we set
\[
E^{\infty}_{p,q} = \varinjlim E^r_{p,q}.
\]

\begin{figure}[ht]
  \begin{center}
    \leavevmode
    \epsfxsize = 10cm
    \epsfysize = 10cm
    \epsfbox{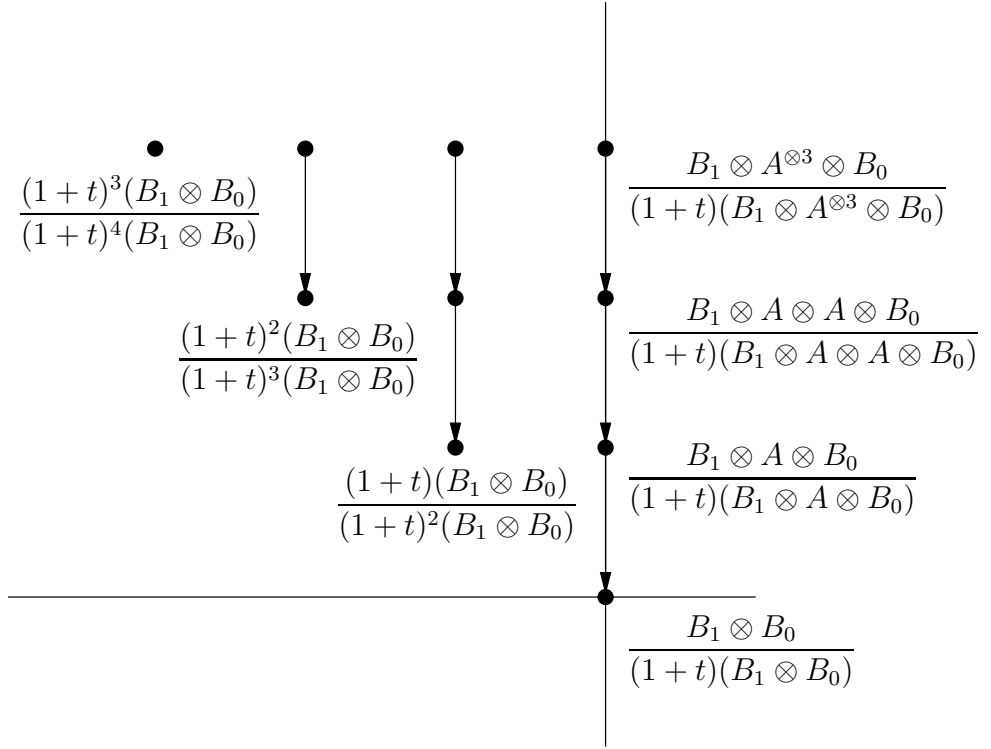}
    \put(-50,35){$\displaystyle{\frac{B_1 \ot B_0}{(1+t)(B_1 \ot B_0)}}$}
    \put(-50,100){$\displaystyle{\frac{B_1 \ot A \ot B_0}{(1+t)(B_1 \ot A \ot B_0)}}$}
    \put(-50,155){$\displaystyle{\frac{B_1 \ot A \ot A \ot B_0}{(1+t)(B_1 \ot A \ot A \ot B_0)}}$}
    \put(-50,210){$\displaystyle{\frac{B_1 \ot A^{\ot 3} \ot B_0}{(1+t)(B_1 \ot A^{\ot 3} \ot B_0)}}$}
    \put(-160,90){$\displaystyle{\frac{(1+t)(B_1 \ot B_0)}{(1+t)^2(B_1 \ot B_0)}}$}
    \put(-220,145){$\displaystyle{\frac{(1+t)^2(B_1 \ot B_0)}{(1+t)^3(B_1 \ot B_0)}}$}
    \put(-280,200){$\displaystyle{\frac{(1+t)^3(B_1 \ot B_0)}{(1+t)^4(B_1 \ot B_0)}}$}
    \caption{The $E^0$ level of the spectral sequence}
    \label{E0-level}
  \end{center}
\end{figure}

The terms at all levels are only nonzero in the second quadrant above
the line $q=-p$.  The $E^0$ level is shown in Figure \ref{E0-level}.
We will be particularly concerned with the terms $E^r_{p,-p}$
along the lower diagonal.

% was [ht]
%\vspace{\li}
%\vspace{-\li}

\begin{lemma}
Every nonzero term of the $E^0$ level of the spectral sequence above
is isomorphic to the tensor product of specialized skein modules.
Namely, for $(p,q)$ with $q \geq 0$ and $-q \leq p \leq 0$, we have
$$E^0_{p,q}
\cong
K_{-1}(H_1) \ot (K_{-1}(F))^{\ot (p+q)} \ot K_{-1}(H_0)$$
\label{E0-lemma}
\end{lemma}

\proof

By definition, the $(p,q)$ term of the $E^0$ level is
\begin{eqnarray}
E^0_{p,q} & = & \mathcal{F}_p \left( C_{p+q} \right) / \mathcal{F}_{p-1} \left( C_{p+q} \right) \nonumber \\
& = & \frac{\displaystyle{\sum_{k=-p}^{\infty} (1+t)^{k}(B_1 \ot A^{\ot (p+q)} \ot B_0)}}
{\displaystyle{\sum_{k=-p+1}^{\infty} (1+t)^{k}(B_1 \ot A^{\ot (p+q)} \ot B_0)}} \nonumber \\
& \cong & \frac{(1+t)^{-p}(B_1 \ot A^{\ot (p+q)} \ot B_0)}{(1+t)^{-p+1}(B_1 \ot A^{\ot (p+q)} \ot B_0)} \nonumber
\end{eqnarray}
and the boundary maps of the $E^0$ level are the Hochschild boundary
maps.  Consider the complex in the column along the $q$-axis in Figure \ref{E0-level}.
These terms are the original Hochschild chains
$C_q = B_1 \ot (A^{\ot q}) \ot B_0$ quotiented by the terms that have a factor
of $(1+t)$ in their coefficients.
Modding out by the ideal generated by $(1+t)$ is the same as setting $(1+t)=0$
or simply evaluating the polynomials at $t=-1$.
Since evaluating at $t=-1$ yields the specialized skein module,
we have, for example, $K(F)/(1+t)K(F) \cong K_{-1}(F)$.
Recall that $A=K(F)$, $B_0=K(H_0)$, and $B_1=K(H_1)$, so
$A/(1+t)A \cong K_{-1}(F)$, $B_0 / (1+t) B_0 \cong K_{-1}(H_0)$,
and $B_1 / (1+t) B_1 \cong K_{-1}(H_1)$.

These quotients are consistent with the tensor product, thus,
\[
\frac{B_1 \ot (A^{\ot n}) \ot B_0}{(1+t) (B_1 \ot (A^{\ot n}) \ot B_0)}
\cong K_{-1}(H_1) \ot (K_{-1}(F))^{\ot n} \ot K_{-1}(H_0).
\]
In addition, for $C_n = B_1 \ot (A^{\ot n}) \ot B_0$ there is a natural map
from $C_n$ to $(1+t)^{-p} C_n$ given by multiplication
by $(1+t)^{-p}$.  Each of the skein modules $A$, $B_0$, and $B_1$
is free on simple diagrams,
so multiplying by $(1+t)^{-p}$ induces an isomorphism
\[ \frac{C_n}{(1+t)C_n} \cong \frac{(1+t)^{-p} C_n}{(1+t)^{-p+1} C_n} \]
and therefore
\begin{eqnarray}
E^0_{p,q} & \cong &
\frac{(1+t)^{-p}(B_1 \ot A^{\ot (p+q)} \ot B_0)}{(1+t)^{-p+1}(B_1 \ot A^{\ot (p+q)} \ot B_0)} \nonumber \\
& \cong & \frac{(B_1 \ot A^{\ot (p+q)} \ot B_0)}{(1+t)(B_1 \ot A^{\ot (p+q)} \ot B_0)} \nonumber \\
& \cong & K_{-1}(H_1) \ot (K_{-1}(F))^{\ot (p+q)} \ot K_{-1}(H_0) \nonumber
\end{eqnarray}
as desired.
\qed

To move from the $E^0$ level to the $E^1$ level in the spectral sequence,
we take the homology of the vertical $E^0$ level complexes.  Since these complexes
are Hochschild complexes, we get the Hochschild homology.  From the proof
of Lemma \ref{E0-lemma}, we know that
\[ \frac{C_n}{(1+t)C_n} \cong \frac{(1+t)^{-p} C_n}{(1+t)^{-p+1} C_n}. \]
It is enough, then, to focus our attention on the complex along the $q$-axis
in Figure \ref{E0-level}.  Namely,
\[
E^1_{0,q} = HH_{q} \left( \frac{A}{(1+t)A};
\frac{B_0}{(1+t) B_0} \ot \frac{B_1}{(1+t)B_1}
\right)
\]
and in general
\begin{equation}
E^1_{p,q} = HH_{p+q} \left( \frac{A}{(1+t)A};
\frac{B_0}{(1+t) B_0} \ot \frac{B_1}{(1+t)B_1}
\right).
\label{e1hom}
\end{equation}

\begin{figure}[ht]
  \begin{center}
    \leavevmode
    \epsfxsize = 10cm
    \epsfysize = 10cm
    \epsfbox{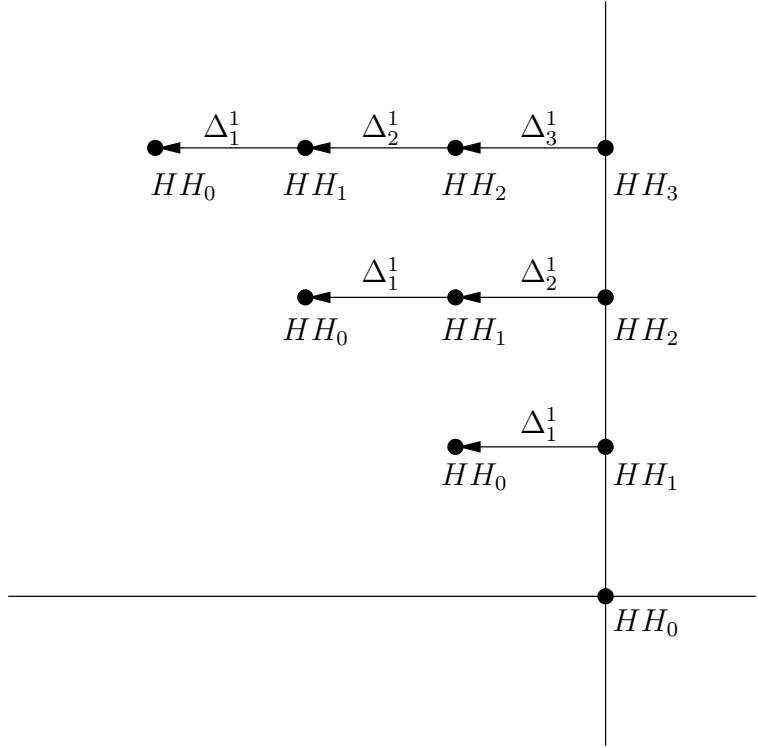}
    \put(-55,45){$HH_0$}
    \put(-55,100){$HH_1$}
    \put(-55,155){$HH_2$}
    \put(-55,210){$HH_3$}
    \put(-90, 120){$\Delta^1_1$}
    \put(-90, 177){$\Delta^1_2$}
    \put(-90, 232){$\Delta^1_3$}
    \put(-120,100){$HH_0$}
    \put(-120,155){$HH_1$}
    \put(-120,210){$HH_2$}
    \put(-150, 177){$\Delta^1_1$}
    \put(-150, 232){$\Delta^1_2$}
    \put(-180,155){$HH_0$}
    \put(-180,210){$HH_1$}
    \put(-210, 232){$\Delta^1_1$}
    \put(-230,210){$HH_0$}
    \caption{The $E^1$ level of the spectral sequence}
    \label{E1-level}
  \end{center}
\end{figure}

The boundary
maps for the $E^1$ level are the induced connecting homomorphisms
$\Delta^1_n : HH_n \to HH_{n-1}$.  The $E^1$ level is shown in
Figure \ref{E1-level}.

%\vspace{\li}
%\vspace{-\li}

The connecting homomorphisms in this context are essentially the same
as the boundary map, we are just considering the image in a different
quotient space.  For example, a $1$-cycle in $E^0_{0,1}$ is an element
$\alpha \in (B_1 \ot A \ot B_0) / (1+t) (B_1 \ot A \ot B_0)$ that maps to
zero in $(B_1 \ot B_0) / (1+t)(B_1 \ot B_0)$, so $d_1(\alpha)$ is
divisible by $(1+t)$.  The connecting homomorphism $\Delta$ is defined
by $\Delta(\alpha) = d_1(\alpha)$, then consider $d_1(\alpha)$
in $(1+t)(B_1 \ot B_0) / (1+t)^2 (B_1 \ot B_0)$.
For example, $\Delta(\alpha)$ will be zero if $d_1(\alpha)$ is
also divisible by $(1+t)^2$.

Thus we have the following lemma.

\begin{lemma}
$\Delta^1_1 : HH_1 \to HH_0$ will be the zero map if every element
whose boundary is divisible by $(1+t)$ also has its boundary divisible
by $(1+t)^2$.  \qed
\end{lemma}

When we move to the next level, we will be looking at something
whose boundary is divisible by $(1+t)^2$ and it will be zero under
the next connecting homomorphism if its boundary is also divisible
by $(1+t)^3$.

Thus, all the connecting homomorphisms will be the zero map if every
element whose boundary is divisible by $(1+t)$ also has a boundary
that is divisible by $(1+t)^r$ for all $r$.  This will happen if
the boundary in question is \textit{exactly} zero, not just some
polynomial divisible by $(1+t)$.

Again, consider Figure \ref{E1-level} and 
notice that these $\Delta$ maps are
horizontal and pointing along the negative $x$ axis.
Recall from Lemma \ref{E0-lemma} that the filtered complex (the $E^0$ level)
can be seen as a complex of specialized skein modules.  Thus the
$E^1_{p,-p}$ terms are just the Hochschild homology of these
specialized skein modules.  Namely,
\begin{eqnarray}
E^1_{p,-p}
& = & HH_0 \left( \frac{A}{(1+t)A} ;
\frac{(B_0 \ot B_1)}{(1+t) (B_0 \ot B_1)} \right) \label{eq-3} \\
& = & HH_0 \left( K_{-1}(F) ; K_{-1}(H_0) \ot K_{-1}(H_1) \right) \nonumber \\
& = & K_{-1}(H_1) \ot_{K_{-1}(F)} K_{-1}(H_0) \nonumber \\
& = & K_{-1}(M). \nonumber
\end{eqnarray}

\section{Looking for torsion}
\label{torsion}

We have constructed a spectral sequence from the
filtered Hochschild complex.  We now want to use this
construction to look for $(1+t)$-torsion in the skein module $K(M)$.

\begin{definition}
A module $X$ over a ring $R$ has \textit{torsion} if there
exist nonzero elements $r \in R$ and $x \in X$ such that
$rx = 0$.
\end{definition}

\begin{definition}
Let $M$ be a manifold. Recall that $K(M)$, the skein module
of $M$, is a module over the Laurent polynomials
$R=\mathbb{Z}[t,t^{-1}]$.
The filtration $\mathcal{F}_s(R) = (1+t)^{-s} R$ for $s \leq 0$
extends to $K(M)$ as
\[
\dots \subset (1+t)^3 K(M) \subset (1+t)^2 K(M) \subset (1+t) K(M)
\subset K(M).
\]
The quotients $K(M) / (1+t)^n K(M)$ together with the projections
$\theta_n : K(M) / (1+t)^n K(M) \to K(M) / (1+t)^{n-1} K(M)$ form an
injective system.
The \textit{completion} of $K(M)$ is the inverse limit of this system,
\[
\overline{K(M)} = \varprojlim K(M) / (1+t)^n K(M).
\]
The completion can also be seen as the module of all sequences
$\{ a_n \}_{n=0}^{\infty}$ with $a_n \in K(M) / (1+t)^n K(M)$ and $\theta_n(a_n) = a_{n-1}$.
There is a homomorphism $\varphi : K(M) \to \overline{K(M)}$ where
$\varphi ( \alpha)$ is the constant sequence $\{ \alpha \}_{n=0}^{\infty}$.
Note that the kernel of this homomorphism is
\[
\mathrm{Ker}(\varphi) = \bigcap_{n=0}^{\infty}(1+t)^n K(M).
\]
\end{definition}

For more information on the definition above and properties of the
completion, see Atiyah and MacDonald \cite[Chapter 10]{AtiyahMacDonald}.
\begin{theorem}
If the $\Delta^r_1: E^r_{p+1,-p} \to E^r_{p,-p}$ maps are identically zero
at every level in the spectral sequence, then there is
no $(1+t)$-torsion in $\overline{K(M)}$.
\label{thm-zero-maps}
\end{theorem}

\proof
We focus our attention in the spectral sequence to the terms
that lie along the lower diagonal (the line $q=-p$).  These
are the terms $E^r_{p,-p}$ for $p \leq 0$.

At the $E^1$ level, from Equation \ref{e1hom} we know that
these terms have the form
\[
E^1_{p,-p} = HH_0 \left( \frac{A}{(1+t) A} ;
\frac{ B_0}{(1+t) B_0} \ot
\frac{ B_1}{(1+t) B_1} \right)
\]
so these terms are the zeroth Hochschild homology modules
of the various filtered complexes.

The maps $\Delta^1_1 : E^1_{p+1,-p} \to E^1_{p,-p}$ are zero
maps, thus
\[E^2_{p,-p} =
\frac{\mathrm{ker}(E^1_{p,-p} \to 0)} {\mathrm{im}(\Delta^1)} = 
E^1_{p,-p}.
\]
That is, the term at position $(p,-p)$ remains unchanged
when we move from the $E^1$ level to the $E^2$ level.

\begin{figure}[ht]
  \begin{center}
    \leavevmode
    \epsfxsize = 10cm
    \epsfysize = 10cm
    \epsfbox{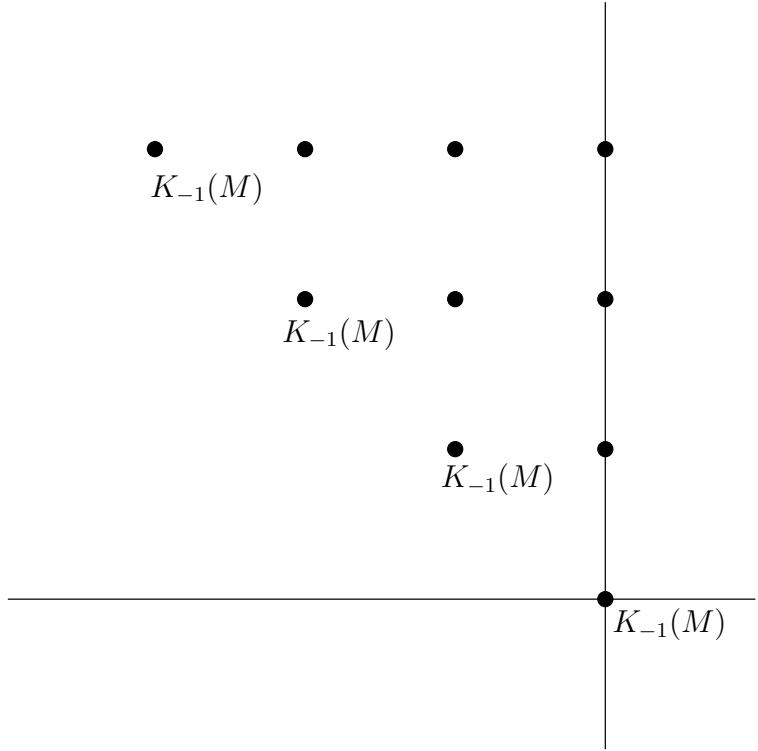}
    \put(-55,45){$K_{-1}(M)$}
    \put(-120,100){$K_{-1}(M)$}
    \put(-180,155){$K_{-1}(M)$}
    \put(-230,210){$K_{-1}(M)$}
    \caption{The $E^{\infty}$ level of the spectral sequence}
    \label{fig-e-infty}
  \end{center}
\end{figure}

The argument is the same for any $r$.
The maps $\Delta^r_1 : E^r_{p+r,-p-r+1} \to E^r_{p,-p}$ are zero
maps, thus
\[E^{r+1}_{p,-p} =
\frac{\mathrm{ker}(E^r_{p,-p} \to 0)} {\mathrm{im}(\Delta^r)} = 
E^r_{p,-p}.
\]

Thus after the $E^1$ level, the terms $E^r_{p,-p}$ along the lower
diagonal are always just the zeroth Hochschild homology of the
filtered complexes, which, from Equation \ref{eq-3}, we know to
be isomorphic to the specialized skein module of $M$.  Thus, in the limit we have
\[
E^{\infty}_{p,-p} = \varinjlim E^r_{p,-p} = E^1_{p,-p} \cong K_{-1}(M).
\]

The $E^{\infty}$ level of the spectral sequence is shown in Figure \ref{fig-e-infty}.
Each term along the lower diagonal is isomorphic to $K_{-1}(M)$.
Specifically,
\[
E^{\infty}_{p,-p} \cong \frac{(1+t)^{-p} K(M)}{(1+t)^{-p+1}K(M)},
\]
and the map
\[ \mu_n : \frac{K(M)}{(1+t)K(M)} \to \frac{(1+t)^nK(M)}{(1+t)^{n+1}K(M)} \]
given as multiplication by $(1+t)^n$ is an isomorphism for all $n \geq 1$.

Suppose there exists an $\alpha \in K(M)$ such that $(1+t) \alpha = 0$.
Then $\mu_1(\alpha) = 0$.
Since $\mu_1$ is an isomorphism, $\alpha$ must be zero in $K(M) / (1+t) K(M)$.
Hence $\alpha$ is divisible by $(1+t)$.
Let $\alpha_1 = \alpha / (1+t)$.  Then
$\mu_2(\alpha_1) = (1+t)^2 \alpha = 0$.  Since $\mu_2$ is an isomorphism, $\alpha_1$
must be zero in $K(M) / (1+t) K(M)$.  Hence $\alpha_1$ is divisible by $(1+t)$ and so
$\alpha$ is divisible by $(1+t)^2$.  An induction argument shows that
$\alpha$ is divisible by $(1+t)^n$ for all $n$.  Thus
$\alpha \in \cap_{n=0}^{\infty} (1+t)^n K(M)$
and so $\alpha = 0$ in the completion $\overline{K(M)}$.
Therefore there is no $(1+t)$-torsion in $\overline{K(M)}$.

%Lastly, we know that the absence of $(1+t)$ torsion in
%$\overline{K(M)}$ implies the absence of torsion in $\overline{K(M)}$.  This
%follows from the fact that the specialized skein module
%$K_{-1}(M)$ is a vector space, hence is torsion free.
%If $f(t)*s = 0$ for a nonzero skein $s$ and a 
%polynomial $f(t)$ nonzero and not divisible by $(1+t)$ we
%have $f(-1)*s = 0$ for a nonzero number $f(-1)$ and a
%nonzero skein $s$.  However, we cannot have such torsion
%in $K_{-1}(M)$ because the specialized skein module is
%a vector space over $\bbc$.
\qed

The motivation for starting the search for torsion by looking
for $(1+t)$-torsion stems from the result of Hoste and Przytycki
\cite{HostePrzytycki1995} about the skein module of
$S^1 \times S^2$.  Namely,
they show that
$$K(S^1 \times S^2) \cong
R \oplus \bigoplus_{i=1}^{\infty} R / (1 - t^{2i+4})$$
where $R = \mathbb{Z}[t, t^{-1}]$.  Since $t=-1$ is a root
of each of the $(1 - t^{2i+4})$ polynomials, all of the torsion in
$K(S^1 \times S^2)$ can be interpreted as $(1+t)$-torsion.  It is
natural to ask if $(1+t)$-torsion is the only kind of torsion that a
skein module can have.  We conclude this section with
the conjecture that $(1+t)$-torsion is indeed
the only kind of torsion.  If true, this conjecture would make
the conclusions of Theorem \ref{thm-zero-maps} more general.

\begin{conjecture}
The absence of $(1+t)$-torsion in $\overline{K(M)}$
implies the absence of torsion in $K(M)$.
\end{conjecture}

\section{Examples}

The existence of torsion in $K(S^1 \times S^2)$ was shown by Hoste and Przytycki
in \cite{HostePrzytycki1995}.  The absence of torsion in the skein module of
each of the other lens spaces was
shown by them in \cite{HostePrzytycki1993}.  Below we show an alternate way to see
$(1+t)$-torsion in $K(S^1 \times S^2)$ and the absence of $(1+t)$-torsion in
$K(L(2,1))$.

\subsection{The space $S^1 \times S^2$.}

Take the genus one Heegaard splitting
$H_1 \cup_{f_1} T^2 \times I \cup_{f_0} H_0$ for $S^1 \times S^2$ where
$f_0 : \partial H_0 \to T^2 \times \{0\}$ is the identity
map and
$f_1 : \partial H_1 \to T^2 \times \{1\}$ is
$$f_1 = \begin{pmatrix} -1 & 0 \\ ~0 & 1 \end{pmatrix}.$$
%$$f_1 = \binom{-1~~0}{~~0~~1}.$$

Consider the element $\alpha \in
\mathrm{Tor}_1^{K_{-1}(T^2)} \Big( K_{-1}(H_1), K_{-1}(H_0) \Big)$ given by
\[
\alpha = \phi \ot (p(t)*y - q(t)*z) \ot \phi = \phi \ot \Big( p(t)*\torusL - q(t)*\torusLM \Big) \ot \phi
\]
with $p(-1) = q(-1) = 1$.
Now consider $\Delta^1_1(\alpha)$ which is given by
\begin{eqnarray}
\Delta^1_1(\alpha) & = & \Big( p(t)*\torusL - q(t)*\torusLMinv \Big) \ot
\phi \nonumber \\
& & - \phi \ot \Big( p(t)*\torusL - q(t)*\torusLM \Big) \nonumber \\
& = & (p(t)+t^3q(t))~\torusL \ot \phi~~~~-~~~~\phi \ot (p(t)+t^{-3}q(t))~\torusL. \nonumber
\end{eqnarray}
Since both $(p(t)+t^3q(t))$ and $(p(t)+t^{-3}q(t))$ are divisible by $(1+t)$, the
element $\alpha$ is a cycle and thus represents a nontrivial element in
$\mathrm{Tor}_1$.  However, $(p(t)+t^3q(t))$ and $(p(t)+t^{-3}q(t))$ are
not divisible by $(1+t)^2$, so $\Delta^1_1(\alpha) \neq 0$.

Now Theorem \ref{thm-zero-maps} does not apply, but if we mimic the proof
of Theorem \ref{thm-zero-maps} for this example
we will see that at the $E^2$ level of the spectral sequence the term
$E^2_{-1,1}$ is no longer isomorphic to $E^1_{-1,1}$.  Instead it is the nontrivial
quotient of $E^1_{-1,1}$ by the image of $E^1_{0,1}$ under $\Delta^1_1$.  Compare
Figure \ref{E1-level}.  Thus at the $E^{\infty}$ level, the map
\[
\frac{K(S^1 \times S^2)}{(1+t)K(S^1 \times S^2)} \to
\frac{(1+t)K(S^1 \times S^2)}{(1+t)^2K(S^1 \times S^2)}
\]
has nontrivial kernel.  Thus there exists an element $\beta \in K(S^1 \times S^2)$
that is not divisible by $(1+t)$ such that $(1+t) \beta$ is divisible by $(1+t)^2$.
Thus $(1+t) \beta = 0$ and $(1+t)$-torsion exists in $K(S^1 \times S^2)$.

\subsection{The $L(2,1)$ lens space.}

Take the genus one Heegaard splitting
$H_1 \cup_{f_1} T^2 \times I \cup_{f_0} H_0$ for $L(2,1)$ where
$f_0 : \partial H_0 \to T^2 \times \{0\}$ is the identity
map and
$f_1 : \partial H_1 \to T^2 \times \{1\}$ is
$$f_1 = \begin{pmatrix} 1 & 2 \\ 1 & 1 \end{pmatrix}.$$
%$$f_1 = \binom{1~~2}{1~~1}.$$
%text below cut from thesis section on S^3, rework as needed to make it
%blend nicely
Let $\ell_i$ and $m_i$ be the longitude and meridian of $H_i$.
Let $\ell$ and $m$ be the longitude and meridian of $T^2 \times I$.

The specialized skein modules of the handlebodies,
$K_{-1}(H_0)$ and $K_{-1}(H_1)$, correspond to subvarieties of
the specialized skein module of the surface,
$K_{-1}(T^2)$.  When we push a framed link from $T^2 \times I$ into
one of the handlebodies, there are relations induced by the fact that
$m_i$ is trivial and $\ell_i \simeq \ell_i m_i$ in handlebody $H_i$.

Recall from Example \ref{ex-torus} that $K_{-1}(T^2)$ is
generated by $x = -\tr(m)$, $y = -\tr(\ell)$, and $z = -\tr(\ell m)$.
Indeed,
$K_{-1}(T^2) = \bbc[x,y,z] / I$ where $I$ is the ideal
$I = (x^2+y^2+z^2+xyz-4)$.
Since $f_0$ is the identity map,
the relations induced by $m_0 \simeq *$ and $\ell_0 \simeq \ell_0 m_0$
are $x = -\tr(m) = -\tr(m_0) = -2$ and
$y = -\tr(\ell) = -\tr(\ell_0) = -\tr(\ell_0 m_0) = -\tr(\ell m) = z.$
Let $J = (x+2, y-z)$, then $K_{-1}(H_0) = K_{-1}(T^2) / J$.

\begin{remark}
In the preceding paragraph (and in the paragraphs below), we are abusing
notation when we use $\tr(m)$, $\tr(\ell)$, etc.  We are suppressing
the representations $\rho : \pi_1(H_i) \to \slc$ and
$\hat{\rho} : \pi_1(T^2) \to \slc$.  The relations on the traces of the
matrices come from the curves themselves.  Thus it seems more instructive
to emphasize the curves over the matrices.
%The subscripts on $m$, $\ell$,
%etc.\ will be used to keep track of the location of the curve.  The curves
%$m_i$, $\ell_i$, and $\ell_i m_i$, for example, are in handlebody $H_i$
%while the curves $m$, $\ell$, and $\ell m$ are in $T^2 \times I$.
\end{remark}

%end of text cut from thesis

The inverse of the gluing map for $H_1$ is
\[ f^{-1}_1 = \begin{pmatrix} -1 & ~2 \\ ~1 & -1 \end{pmatrix}.\]
%\[f^{-1}_1 = \binom{-1~~~~~~2}{~~~1~~-1}.\]
Thus the inclusion of $T^2$ into $H_1$
sends $\ell^2 m$ to $m_1$, $\ell$ to $\ell_1^{-1} m_1$, and
$\ell m$ to $\ell_1$.  The relation induced by $m_1 \simeq *$ uses $f_1^{-1}$.
In particular, $m_1 \simeq *$ induces
$\tr(\ell^2 m) = \tr(m_1) = 2$.  Using the trace identity for
$\slc$ we have
$\tr(\ell^2 m) = \tr(\ell)*\tr(\ell m) - \tr(m)$.
Thus $\tr(\ell^2 m) = 2$ becomes $yz+x = 2.$  Similarly a relation is induced
by $\ell_1 \simeq \ell_1 m_1 \simeq \ell_1^{-1} m_1$ using $f_1^{-1}$.
In particular,
\[y = -\tr(\ell) = -\tr(\ell_1^{-1} m_1) = -\tr(\ell_1) = -\tr(\ell m) = z.\]
Let $K = (yz+x-2, y-z)$, then $K_{-1}(H_1) = K_{-1}(T^2) / K$.

\begin{lemma}
As a vector space over $\bbc$, $\mathrm{Tor}_1^{K_{-1}(T^2)} \Big( K_{-1}(H_1),
K_{-1}(H_0) \Big)$ is spanned by the set $\{ (y-z), y(y-z) \}$.
\end{lemma}

\proof
We know that $\mathrm{Tor}_1^{K_{-1}(T^2)} \Big( K_{-1}(H_1),
K_{-1}(H_0) \Big) = (J \cap K)/(JK)$ with $J = (x+2, y-z)$ and
$K = (yz+x-2, y-z)$.  Take $\alpha(x,y,z) \in (J \cap K)$.
We know that \[ \alpha(x,y,z) = p_1(x,y,z) (x+2) + p_2(x,y,z) (y-z) \] and
\[ \alpha(x,y,z) = q_1(x,y,z)(yz+x-2) + q_2(x,y,z) (y-z). \]
Use $(x+2)(yz+x-2) \in JK$ and $(y-z)(yz+x-2) \in JK$ to write $q_1(x,y,z)$
as a function in $y$.    Use $(x+2)(y-z)$ and $(y-z)(y-z)$ to write
$q_2(x,y,z)$ as a function in $y$.  Thus
\[ \alpha(x,y,z) = \tilde{q}_1(y)(yz+x-2) + \tilde{q}_2(y)(y-z) \]
in the quotient $(J \cap K)/(JK)$.
Evaluating $\alpha(x,y,z)$ at $(-2,y,y)$ we have
\[0 = \alpha(-2,y,y) = \tilde{q}_1(y)(y^2-4) + \tilde{q}_2(y)(y-y) = \tilde{q}_1(y)(y^2-4).\]
Thus $\tilde{q}_1(y) = 0$ and $\alpha(x,y,z) = \tilde{q}_2(y)(y-z)$.  Hence
$(J \cap K)/(JK)$ is spanned by the set $\{ (y-z), y(y-z), y^2(y-z), y^3(y-z), \dots \}$.

Now consider the element $y^2(y-z) - 4(y-z)$ in $(J \cap K)/(JK)$ as follows.
\begin{eqnarray}
y^2(y-z) - 4(y-z) & = & y^2(y-z) - 2(y-z) - 2(y-z) + (x+2)(y-z)\nonumber \\
& = & y^2(y-z) - y(y-z)(y-z) - 2(y-z) + x(y-z)\nonumber \\
& = & y^2(y-z) - y^2(y-z) + yz(y-z) - 2(y-z) + x(y-z)\nonumber \\
& = & (yz + x - 2)(y-z).\nonumber
\end{eqnarray}
Since $(yz+x-2)(y-z) \in JK$, we have $y^2(y-z) = 4(y-z)$ in $(J \cap K)/(JK)$.
Thus $(J \cap K)/(JK)$ is spanned by the set $\{ (y-z), y(y-z) \}$.
\qed

\begin{proposition}
There is no $(1+t)$-torsion in $\overline{K(L(2,1))}$.
\end{proposition}

\proof
%As in the proof of Theorem \ref{thm-torsion-s3}
Any element
in $(J \cap K)/(JK)$ can be written as $\phi \ot \alpha \ot \phi$
where $\alpha$ is a linear combination of $(y-z)$ and $y(y-z)$.  To show
that $\Delta_1^r (\phi \ot \alpha \ot \phi) = 0$ it is enough to show that
$\Delta_1^r (\phi \ot (y-z) \ot \phi) = 0$ and $\Delta_1^r (\phi \ot y(y-z) \ot \phi) = 0$.

In $A / (1+t) A$ the element $y-z$ is equal to $y+t^3z$.
\begin{eqnarray}
\Delta_1^r ( \phi \ot (y + t^3 z) \ot \phi ) & = &
f_1^{-1} ( \ell + t^3 \ell m) \ot \phi - \phi \ot (\ell_0 + t^3 \ell_0 m_0) \nonumber \\
& = & (\ell_1^{-1} m_1 + t^3 \ell_1) \ot \phi - \phi \ot (\ell_0 + t^3 \ell_0 m_0 ) \nonumber \\
& = & (-t^3 \ell_1 + t^3 \ell_1) \ot \phi - \phi \ot (\ell_0 +(t^3)(-t^{-3}) \ell_0) \nonumber \\
& = & 0 \ot \phi - \phi \ot 0 \nonumber \\
& = & 0 \nonumber
\end{eqnarray}

In $A / (1+t) A$, the element $y(y-z)$ is equal to
$\beta = p(t)*y^2 - \frac{1}{2} yz - \frac{1}{2} zy$
where $p(t) = - \frac{1}{2} t^{-3} - \frac{1}{2} t^{-5}$.
Apply $\Delta^r_1$ to the element $\phi \ot \beta \ot \phi$.

% i'm not sure why, but the equation below seems to have
% two extra spaces above it.  if this problem persists, it can
% be corrected with the folowing vspace line
%\vspace{-2\li}
\begin{eqnarray}
\Delta^r_1 ( \phi \ot \beta \ot \phi) & = &
\Delta^r_1 \Big( \phi \ot \Big( p(t) y^2 - {\textstyle \frac{1}{2}} yz - {\textstyle \frac{1}{2}} zy \Big) \ot \phi \Big) \nonumber \\
& = & f_1^{-1} \Big( p(t) \ell*\ell - {\textstyle \frac{1}{2}} \ell * \ell m - {\textstyle \frac{1}{2}} \ell m * \ell \Big) \ot \phi \nonumber \\
& & - \phi \ot \Big( p(t) \ell_0 * \ell_0 - {\textstyle \frac{1}{2}} \ell_0 * \ell_0 m_0 - {\textstyle \frac{1}{2}} \ell_0 m_0 * \ell_0 \Big) \nonumber \\
& = & \Big( p(t) \ell_1 m_1^{-1} * \ell m_1^{-1} - {\textstyle \frac{1}{2}} \ell_1 * \ell_1 m_1^{-1}
- {\textstyle \frac{1}{2}} \ell_1 m_1^{-1} * \ell_1 \Big) \ot \phi \nonumber \\
& & - \phi \ot \Big( p(t)*\ell_0 * \ell_0 - {\textstyle \frac{1}{2}} \ell_0 * \ell_0 m_0 - {\textstyle \frac{1}{2}} \ell_0 m_0 * \ell_0 \Big)
\label{eq-56}
\end{eqnarray}

\begin{figure}[ht]
  \begin{center}
    \leavevmode
    \epsfxsize = 4cm
    \epsfysize = 2cm
    \epsfbox{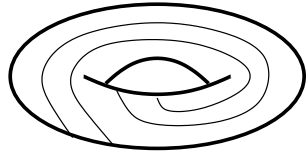}
    \caption{The link $\gamma$ in a solid torus}
    \label{fig-gamma}
  \end{center}
\end{figure}

\begin{figure}[ht]
  \begin{center}
    \leavevmode
    \epsfxsize = 4cm
    \epsfysize = 2cm
    \epsfbox{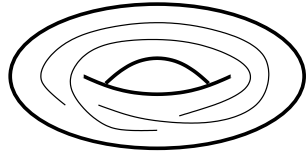}
    \caption{The link $\delta$ in a solid torus}
    \label{fig-delta}
  \end{center}
\end{figure}

\begin{figure}[ht]
  \begin{center}
    \leavevmode
    \epsfxsize = 12cm
    \epsfysize = 5cm
    \epsfbox{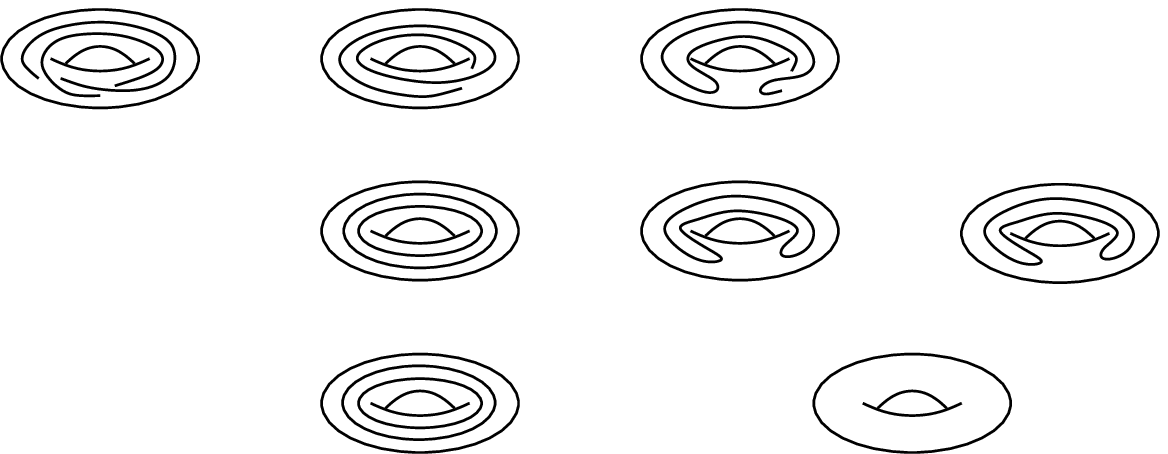}
    \put(-280, 123){$=~~~t$}
    \put(-185, 123){$+~t^{-1}$}
    \put(-280, 68){$=~~t^2$}
    \put(-185, 68){$~+$}
    \put(-90,68){$-~t^{-4}$}
    \put(-280,11){$=~~t^2$}
    \put(-185,11){$+~~(t^{-4} - 1) [2]$}
    \caption{The skein relations applied to the link $\delta$}
    \label{fig-skein-delta}
  \end{center}
\end{figure}

Let $\gamma = \ell_1 m_1^{-1} * \ell_1 m_1^{-1}$ as shown in Figure \ref{fig-gamma}.
Removing the kinks in $\gamma$, we see that $\gamma = t^6 \delta$ where $\delta$
is the link shown in Figure \ref{fig-delta}.  Note also that $\ell_1 m_1^{-1} * \ell_1 = -t^3 \delta$
and $\ell_0 m_0 * \ell_0 = -t^{-3} \bar{\delta}$ where $\bar{\delta}$ is the mirror image of
$\delta$.
%\vspace{\li}
%\vspace{-\li}
%\vspace{\li}
%\vspace{-\li}
Removing kinks and using $\gamma$ and $\delta$, Equation \ref{eq-56} becomes
\begin{eqnarray}
\Delta^r_1( \phi \ot \beta \ot \phi ) & = &
\Big( p(t) t^6 \delta + \LittleOneHalf t^3 \ell_1 * \ell_1 + \LittleOneHalf t^3 \delta \Big) \ot \phi \nonumber \\
& & - \phi \ot \Big( p(t) \ell_0 * \ell_0 + \LittleOneHalf t^{-3} \ell_0 * \ell_0 + \LittleOneHalf t^{-3} \bar{\delta} \Big)
\label{eq-57}
\end{eqnarray}
Now apply the skein relations to $\delta$ as shown in Figure \ref{fig-skein-delta}.
In handlebody $H_1$ we have
$\delta = t^2 \ell_1 * \ell_1 + (t^{-4} - 1)[2] \phi$ and in handlebody $H_0$ we have
$\bar{\delta} = t^{-2} \ell_0 * \ell_0 + (t^4 - 1)[2] \phi$.

Let
$ \beta_1 = p(t) t^6 \delta + \LittleOneHalf t^3 \ell_1 * \ell_1 + \LittleOneHalf t^3 \delta$
and
$ \beta_0 = p(t) \ell_0 * \ell_0 + \LittleOneHalf t^{-3} \ell_0 * \ell_0 + \LittleOneHalf t^{-3} \delta.$
Then $\Delta( \phi \ot \beta \ot \phi ) = \beta_1 \ot \phi - \phi \ot \beta_0$.
Consider
\begin{eqnarray}
\beta_1 & = & p(t) \Big( t^8 \ell_1 * \ell_1 + t^6(t^{-4} - 1)[2] \phi \Big) +
\LittleOneHalf t^3 \ell_1 * \ell_1 + \LittleOneHalf t^5 \ell_1 * \ell_1 +
\LittleOneHalf t^3(t^{-4} - 1)[2] \phi \nonumber \\
& = & (-\LittleOneHalf t^{-3} - \LittleOneHalf t^{-5})t^6(t^{-4} - 1)[2] \phi
+ \LittleOneHalf t^3(t^{-4} - 1)[2] \phi \nonumber \\
& = & - \LittleOneHalf t (t^{-4} - 1)[2] \phi \nonumber
\end{eqnarray}
and
\begin{eqnarray}
\beta_0 & = & p(t) \ell_0 * \ell_0 +
\LittleOneHalf t^{-3} \ell_0 * \ell_0 + \LittleOneHalf t^{-3} \Big( t^{-2} \ell_0 * \ell_0 +
(t^{4} - 1)[2] \phi \Big) \nonumber \\
& = & \LittleOneHalf t^{-3}(t^{4} - 1)[2] \phi. \nonumber
\end{eqnarray}
Then
\begin{eqnarray}
\Delta^r_1( \phi \ot \beta \ot \phi ) & = & \beta_1 \ot \phi - \phi \ot \beta_0 \nonumber \\
& = & -\LittleOneHalf t (t^{-4} - 1)[2] \phi \ot \phi - \phi \ot \LittleOneHalf t^{-3}(t^{4} -1)[2] \phi
\nonumber \\
& = & (- \LittleOneHalf t^{-3} + \LittleOneHalf t)[2] \phi \ot \phi - \phi \ot (\LittleOneHalf t
- \LittleOneHalf t^{-3}) [2] \phi \nonumber \\
& = & 0. \nonumber
\end{eqnarray}

Therefore the $\Delta^r_1$ maps are all zero maps and by
Theorem \ref{thm-zero-maps} there is no $(1+t)$-torsion in
$\overline{K(L(2,1))}$.
\qed

%\subsection{The Poincar\'{e} Manifold}

\section{Torsion in an Homology Sphere}

The computational methods detailed in the previous section are
cumbersome.  Indeed, the description of the skein modules of the lens spaces
given by Hoste and Przytycki in \cite{HostePrzytycki1993, HostePrzytycki1995} is cleaner.
However, the framework of Hochschild homology given by the current paper
will hopefully allow us to use advanced ideas from homological algebra
and representation theory
to search for torsion in the skein module of a manifold.  In
particular, we hope to use the results of Serre in \cite{Serre} and of Goldman and Millson
in \cite{GoldmanMillson1988} to prove the following rather ambitious conjecture.

\begin{conjecture}
If $M$ is an homology sphere, then there is no
$(1+t)$-torsion in $\overline{K(M)}$.
\end{conjecture}

\bibliography{torsion}		% pulls biblio info from the torsion.bib file

\begin{thebibliography}{10}

\bibitem{AtiyahMacDonald}
M.~F. Atiyah and I.~G. MacDonald.
\newblock {\em Introduction to Commutative Algebra}.
\newblock Addison-Wesley, 1969.

\bibitem{Brylinski1988}
J.-L. Brylinski.
\newblock A differential complex for {P}oisson manifolds.
\newblock {\em Journal of Differential Geometry}, 28:93--114, 1988.

\bibitem{Bullock1997}
D.~Bullock.
\newblock Rings of ${SL}_2(\mathbb{C})$-characters and the {K}aufmann bracket
  skein module.
\newblock {\em Commentarii Mathematici Helvetici}, 72:521--542, 1997.

\bibitem{CullerShalen1983}
M.~Culler and P.~B. Shalen.
\newblock Varieties of group representations and splittings of 3-manifolds.
\newblock {\em Annals of Mathematics}, 117:109--146, 1983.

\bibitem{FrohmanGelca2000}
C.~Frohman and R.~Gelca.
\newblock Skein modules and the noncommutative torus.
\newblock {\em Transactions of the American Mathematical Society},
  352(10):4877--4888, 2000.

\bibitem{GoldmanMillson1988}
W.~M. Goldman and J.~J. Millson.
\newblock The deformation theory of representations of fundamental groups of
  compact {K}\"{a}hler manifolds.
\newblock {\em Publications math\'{e}matiques de l'I.H.\'{E}.S.}, 67:43--96,
  1988.

\bibitem{HostePrzytycki1993}
J.~Hoste and J.~H. Przytycki.
\newblock The $(2, \infty)$-skein module of lens spaces; a generalization of
  the {J}ones polynomial.
\newblock {\em Journal of Knot Theory and Its Ramifications}, 2(3):321--333,
  1993.

\bibitem{HostePrzytycki1995}
J.~Hoste and J.~H. Przytycki.
\newblock The {K}auffman bracket skein module of ${S}^1 \times {S}^2$.
\newblock {\em Mathematische Zietschrift}, 220:65--73, 1995.

\bibitem{KauffmanLins}
L.~H. Kauffman and S.~L. Lins.
\newblock {\em Temperley-Lieb Recoupling Theory and Invariants of 3-Manifolds}.
\newblock Princeton University Press, 1994.

\bibitem{Lickorish1993}
W.~B.~R. Lickorish.
\newblock The skein method for 3-manifold invariants.
\newblock {\em Journal of Knot Theory and Its Ramifications}, 2(2):171--194,
  1993.

\bibitem{Osborne}
M.~S. Osborne.
\newblock {\em Basic Homological Algebra}.
\newblock Springer, 2000.

\bibitem{Przytycki1991}
J.~H. Przytycki.
\newblock Skein modules of $3$-manifolds.
\newblock {\em Bulletin of the Polish Academy of Sciences}, 39(1-2):91--100,
  1991.

\bibitem{PrzytyckiSikora2000}
J.~H. Przytycki and A.~Sikora.
\newblock On skein algebras and ${S}l_2(\mathbb{C})$-character varieties.
\newblock {\em Topology}, 39:115--148, 2000.

\bibitem{Rolfsen}
D.~Rolfsen.
\newblock {\em Knots and Links}.
\newblock Publish or Perish, Inc., 1976.

\bibitem{Serre}
J.-P. Serre.
\newblock {\em Local Algebra}.
\newblock Springer, 2000.

\bibitem{Turaev1988}
V.~G. Turaev.
\newblock The {C}onway and {K}auffman modules of the solid torus.
\newblock {\em Zapiski Nauchnykh Seminarov (LOMI)}, 167:79--89, 1988.
\newblock English translation: Journal of Soviet Mathematics, 52(1):2799-2805,
  1990.

\end{thebibliography}
\bibliographystyle{plain}

\end{document}